\documentclass[12pt]{amsart}
\usepackage[left=1in,top=1in,right=1in,bottom=1in,letterpaper]{geometry}
\usepackage{graphics}
\usepackage{amsmath}
\usepackage{amssymb}
\usepackage{amscd}
\usepackage{url}
\usepackage{verbatim} 
\usepackage{color}
\usepackage{stmaryrd}
\usepackage{amsthm}
\usepackage{hyperref}

\title{Another Veech triangle}
\author[]{W. Patrick Hooper}
\thanks{Supported by N.S.F. Postdoctoral Fellowship DMS-0803013}
\address{
The City College of New York\\
New York, NY, USA 10031}
\email{whooper@ccny.cuny.edu}

\newtheorem{theorem}{Theorem}

\newtheorem{lemma}[theorem]{Lemma}

\newtheorem{discussion}[theorem]{Discussion}

\newtheorem{openquestion}[theorem]{Open Question}

\def\endproof{$\diamondsuit$  \newline}

\def\C{\mathbb{C}}% 
\def\R{\mathbb{R}}% 
\def\Z{\mathbb{Z}}% 

\def\GL{\textit{GL}}
\def\SL{\textit{SL}}

\def\H{\mathbb H}%
\def\isomS1{\textit{Isom}_+(S^1)}
\def\rt3{\sqrt{3}}
\def\RAB{R_{\overline{AB}}}
\def\RAE{R_{\overline{AE}}}
\def\RBC{R_{\overline{BC}}}
\def\RCD{R_{\overline{CD}}}
\def\RDE{R_{\overline{DE}}}

\def\MT{\textit{S}}
\def\Y{\textit{Y}}

\def\SLhat{\widehat{\SL}}%

\begin{document}
\bibliographystyle{amsalpha}
\begin{abstract}
We show that the triangle with angles $\frac{\pi}{12}$, $\frac{\pi}{3}$, and $\frac{7 \pi}{12}$ has the lattice property and compute this triangle's Veech group.
\end{abstract}
\maketitle

In a 1989 paper of Veech, a property of a polygon $P$ in the plane was found which implies that the number of combinatorially distinct periodic billiard paths of length less than
$t$, $N(P,t)$, satisfies 
\begin{equation}
\lim_{t \to \infty} \frac{1}{t^2} N(P,t)=c(P),
\end{equation}
where $c(P)$ is a constant depending only on the polygon \cite{V}. This property is now known as the {\em lattice property}. Since that time, the question of which polygons have the lattice property has been studied.
The goal of this paper is to demonstrate that the triangle $\Delta$ with angles
$(\frac{\pi}{12}, \frac{\pi}{3}, \frac{7 \pi}{12})$ has Veech's lattice property. 

\begin{figure}[h]
\begin{center}
\input{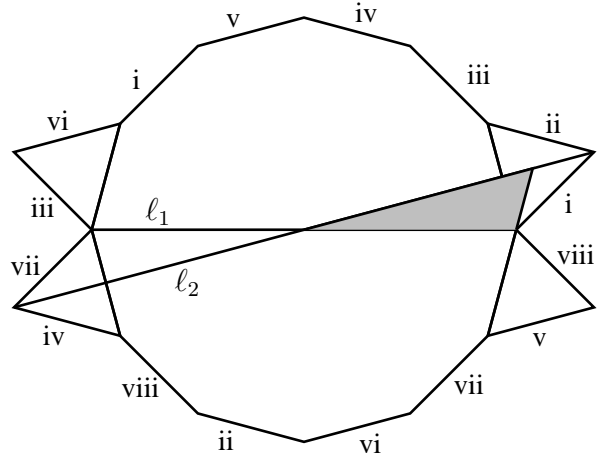}
\caption[The translation surface for $\Delta=(\frac{\pi}{12}, \frac{\pi}{3}, \frac{7 \pi}{12})$]
{The triangle $\Delta$ together with 
the corresponding translation surface, $\MT_\Delta$. Roman numerals indicate edge identifications. }
\label{fig:surface}
\end{center}
\end{figure}

The Zemljakov-Katok construction associates a polygon $P$ with a surface $S_P$ \cite{ZK}. To construct $S_P$, we define a subgroup $G$ of the orthogonal group $O(2)$ generated by reflections in the sides of $P$. Then we define 
$$S_P=\bigsqcup_{g \in G} g(P) / \sim,$$
where $\sim$ is an equivalence relation defined by gluing edges. For every edge $e \subset P$, we glue $g_1(P)$ and $g_2(P)$ by an orientation reversing isometry 
along the edges $g_1(e)$ and $g_2(e)$ whenever $g_1 g_2^{-1}$ is the element of $O(2)$ obtained by reflection in the side $e$ of $P$. These edge identifications glue edges of polygons by translations. Assuming the group $G$ is finite, the surface $S_P$ has finite area. Figure \ref{fig:surface}
shows the triangle $\Delta$ together with the surface $\MT_\Delta$.

The surface $\MT_\Delta$ fits into an infinite class of surfaces with the lattice property of genus $4$ 
discovered by McMullen \cite{McM06}. See the next section for further historical details.

The surface $S_P$ is a {\em translation surface}, a surface built from polygons in $\R^2$ with edges glued by translations. Translation surfaces typically have singular points which are cone points with cone angles that are integer multiples of $2 \pi$. An {\em affine automorphism} of a translation surface $S$ is a homeomorphism $\phi:S \to S$ which preserves the underlying affine structure. By identifying tangent planes of non-singular points with the plane, we see that an affine automorphism 
has a well defined derivative $D(\phi):\R^2 \to \R^2$ and that
\begin{equation}
\label{eq:SLhat}
D(\phi) \in \SLhat(2,\R)=\{M \in \GL(2,\R) | \textit{Det}(M)=\pm 1\}.
\end{equation}

The {\em Veech group} $\Gamma(S) \subset \SLhat(2,\R)$ of a translation surface $S$ is the group of derivatives of affine automorphisms of $S$. The surface $S$ has the {\em lattice property} if $\H^2/\Gamma(S)$ has finite hyperbolic area.

In this paper, we prove the following theorem.

\begin{figure}[htbp]
\begin{center}
\input{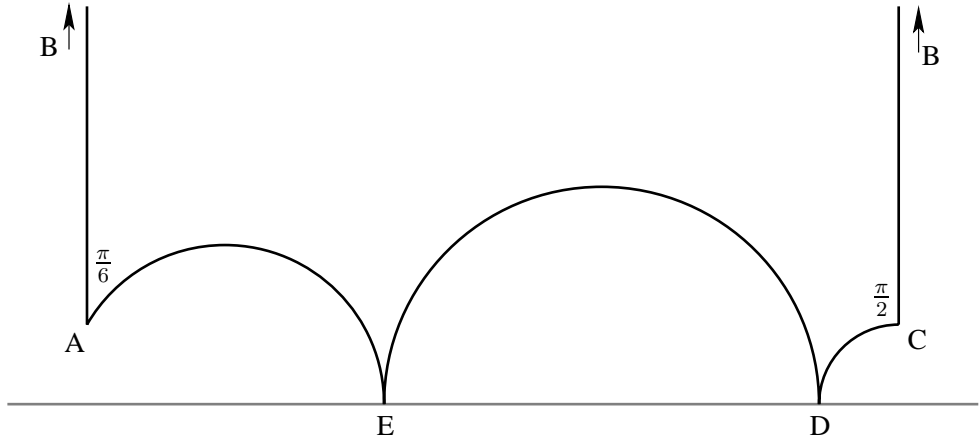}
\caption[The fundamental domain for the action on $\H^2$]{
The fundamental domain for the action of $\Gamma(\MT_\Delta)$ on the
upper half-plane is the polygon pictured. 
The polygon is the hyperbolic convex
hull of its vertices:
$A=i$,
$B=\infty$,
$C=5+3 \rt3 +i$,
$D=4+3 \rt 3$, and
$E=2+\rt3$.
}
\label{fig:domain}
\end{center}
\end{figure}

\begin{theorem}
\label{thm:main}
The surface $\MT_\Delta$ has the lattice property. A fundamental domain for the 
action of the affine automorphism group $\Gamma(\MT_\Delta)\subset \SLhat(2,\R)$ on the hyperbolic plane 
is shown in figure \ref{fig:domain}. $\Gamma(\MT_\Delta)$ is generated by reflections in the sides of this polygon together with $-I$.
\end{theorem}

\section{Historical Remarks}

Theorem \ref{thm:main} adds an additional triangle to the list of known triangles with the lattice property. This list follows. 
\begin{enumerate}
\item The acute isosceles triangles with angles $(\frac{(n-1)\pi}{2n},\frac{(n-1)\pi}{2n},\frac{\pi}{n})$ for $n \ge 3$.
(due to Veech \cite{V}).
\item The acute triangles
$(\frac{\pi}{4}, \frac{\pi}{3}, \frac{5 \pi}{12})$, 
$(\frac{\pi}{5}, \frac{\pi}{3}, \frac{7 \pi}{15})$, 
$(\frac{2 \pi}{9}, \frac{\pi}{3}, \frac{4 \pi}{9})$.
(due to Veech \cite{V}, Vorobets \cite{Vo}, and Kenyon and Smillie \cite{KS} respectively).
\item The right triangles with angles $(\frac{\pi}{n},\frac{(n-2)\pi}{2n},\frac{\pi}{2})$ for $n \ge 4$. (due to Veech \cite{V}).
\item The obtuse isosceles triangles with angles $(\frac{\pi}{n},\frac{\pi}{n},\frac{(n-2)\pi}{n})$ with
$n \ge 5$. (due to Veech \cite{V}).
\item The obtuse triangles with angles $(\frac{\pi}{2n}, \frac{\pi}{n}, \frac{(2n-3)\pi}{2n})$
with $n \ge 4$. (due to Ward \cite{Ward}).
\end{enumerate}

The work of Kenyon and Smillie \cite{KS} together with that of 
Puchta \cite{P} has shown that this list is complete in the case of acute, right, and isosceles triangles.
Kenyon and Smillie provided a simple criterion which can be used to decide that a given triangle does not have
the lattice property. This criterion requires knowledge of the existence of a periodic billiard trajectory. Kenyon and Smillie used the Fagnano curve, which is a periodic billiard path in every acute triangle, to rule out all but a short list of acute triangles with angles that are rational multiples of $\pi$ with denominator less than 10,001. Puchta later eliminated the remaining acute triangles.

The triangle $\Delta$ with angles $(\frac{\pi}{12}, \frac{\pi}{3}, \frac{7 \pi}{12})$ was found using the methods
of Kenyon and Smillie. Rich Schwartz and the author have written a computer program called
McBilliards\footnote{McBilliards is freely available from \url{http://mcbilliards.sourceforge.net/}.}
which is capable of finding periodic billiard trajectories in triangles. This program together with Kenyon and Smillie's
criterion enabled the author to (non-rigorously) search for triangles which might satisfy the lattice property.
The triangle $\Delta$ was the only new triangle with angles of small denominator which seemed to pass this test. (The test did not check for certain
numerical errors.) This suggests that the answer to the following question is likely to be ``yes.''

\begin{openquestion}
Is the list of obtuse triangles with the lattice property now complete?
\end{openquestion}

This triangle was found by the author in 2005 and the result was included in the author's PhD thesis \cite{Hthesis}. As this surface was discovered by the author, McMullen was finishing work on the paper \cite{McM06} which provides infinite lists of surfaces with the lattice property in genera 3 and 4. McMullen was able to check that the surface $\MT_\Delta$ fit into his infinite list. See also the remarks in \cite{McM06} at the end of section 1. 

Thus there are two proofs that $\MT_\Delta$ has the lattice property. In \cite{McM06}, McMullen gives a proof which does not explicitly yield the Veech group but works for a more general class of surfaces.
Our main motivation is to explicitly compute the Veech group. We can see from this computation that the triangle $\Delta$ is special: it is only triangle with the lattice property that has a Veech group which is not a triangle group.

\section{Outline}

In the following section, we will define the affine automorphism group and give enough background to prove the theorem.
The theorem's proof lies in section \ref{sect:proof}. In the final section, we will describe another affinely equivalent translation surface
with Euclidean symmetry group of order $8$. This alternate picture of the surface realizes $\MT_\Delta$ as one
of McMullen's X-shaped polygons. See \cite{McM06}.

\section{Background}

In this section, we briefly define a translation surface, its affine automorphism group, and the lattice property. For more details see
\cite{MT}.

A {\em translation surface} is a closed oriented surface $\MT$ together with a discrete set $\Sigma \subset \MT$ and an atlas of charts 
from $\MT \smallsetminus \Sigma$ to the plane so that the transition functions are translations. The subset $\Sigma$
is a known as the {\em singular set}.
The {\em atlas of charts} is a covering of $\MT \smallsetminus \Sigma$ by open sets
$U_i$ together with local homeomorphisms $\phi_i:U_i \to \R^2$. The {\em transition functions}
are the maps $\phi_i \circ \phi_j^{-1}$ restricted to $\phi_j(U_i \cap U_j)$. 
The translation surface inherits the pull back metric from the plane and also the notion of direction.
Small open sets of $\MT \smallsetminus \Sigma$ are thus isometric to the plane and the points of $\Sigma$ are 
cone points that have cone angles which are integer multiples of $2 \pi$. The relevant example of a translation surface is shown in 
figure \ref{fig:surface}.

We will let $\SLhat(2,\R)$ be the subgroup of affine transformations of the plane that preserve area and fix the origin. See equation \ref{eq:SLhat}.
An element $A \in \SLhat(2,\R)$ acts affinely on the plane. Given $\MT$ we can form a new translation surface $A(\MT)$ by post composing the charts of
$\MT$ with $A$. The transition functions of $A(\MT)$ are translations, since they are just the transition functions of $\MT$ conjugated by $A$. Thus, $A(\MT)$ is another
translation surface.

The {\em affine automorphism group}, $\Gamma(\MT) \subset \SLhat(2,\R)$, of $\MT$ is the set of elements $A \in \SLhat(2,\R)$ so that there is a
direction preserving isometry $\phi:\MT \to A(\MT)$. (Direction preserving is important, otherwise rotations would automatically be in $\Gamma(\MT)$.)
A translation surface $\MT$ is said to have the {\em lattice property} if $\Gamma(\MT) \subset \SLhat(2,\R)$ is a lattice.

Veech discovered a relevant and powerful lemma about parabolics in the affine automorphism group in terms of cylinders of the surface. See \cite{V} and \cite{MT}.
The {\em modulus} of a cylinder is the height of the cylinder divided by its circumference.

\begin{lemma}[Veech]
There is a parabolic in the group $\Gamma(\MT)$ fixing the direction $\theta$ if and only if there is a decomposition of the surface into cylinders in the direction $\theta$ 
whose moduli are commensurable (rational multiples of one another). 
\end{lemma}

The following makes this lemma more explicit.

\begin{discussion}
Suppose $\theta$ is the horizontal direction and $\alpha$ is the greatest common divisor of the moduli of the cylinders in the horizontal decomposition given by the lemma.
The greatest common divisor of a set of commensurable numbers $\{m_1, \ldots, m_n\}$ is the largest number $\alpha$ so that $\frac{m_i}{\alpha} \in \Z$ for all $i$.
The generating parabolic fixing $\theta$ is given by 
$$\left(
\begin{array}{cc}
1 & \frac{1}{\alpha} \\
0 & 1
\end{array}
\right).
$$
\end{discussion}

To aid in visualizing $\Gamma(\MT)$, it is worth considering the action of $\SLhat(2,\R)$ on the upper half plane.
The group $\widehat{\textit{SL}}(2,\R)$ of equation \ref{eq:SLhat} acts on the upper half-plane
by hyperbolic isometries in the standard way. The upper half-plane is a subset of the Riemann sphere, 
$\hat{\C}=\C^2/(\C \smallsetminus \{0\})$. The upper half plane is the equivalence classes of 
elements $(z,1) \in \C^2$ where $z$ has positive imaginary part. An element of
$\widehat{\textit{SL}}(2,\R)$ acts on $\hat{\C}$ as follows:
\begin{equation}
\left(
\begin{array}{cc}
a & b \\
c & d 
\end{array}
\right)
\left(
\begin{array}{c}
z \\
w 
\end{array}
\right)=\left\{
\begin{array}{ll}
\left(
\begin{array}{c}
a z + b w \\
c z + d w 
\end{array}
\right)
& \textit{ if } ad-bc=1 \\
\left(
\begin{array}{c}
a\overline{z} + b\overline{w} \\
c\overline{z} + d\overline{w} 
\end{array}
\right)
& \textit{ if } ad-bc=-1
\end{array}
\right.
\end{equation}
Note that this action is not faithful because $-I$ acts trivially.

\section{The Proof}
\label{sect:proof}
We break up the proof of the theorem into two lemmas. In the first we prove that the elements we list are in $\Gamma(\MT_\Delta)$.
Then we will show that this list generates all of $\Gamma(\MT_\Delta)$.

\begin{lemma}
Each of the reflections in the side of the polygon of figure \ref{fig:domain} is in $\Gamma(\MT_\Delta)$.
$-I$ is also in $\Gamma(\MT_\Delta)$.
\end{lemma}
{\em Proof:} The surface $\MT_\Delta$ has several Euclidean automorphisms. 
The element $-I$ acts on the plane by a Euclidean rotation by $\pi$. Thus,
it is clear $-I \in \Gamma(\MT_\Delta)$.
The Euclidean automorphism group is generated by reflections in lines 
$\ell_1$ and $\ell_2$ of figure \ref{fig:surface}. 
This gives two of our generators:
\begin{equation}
\RAB=\left(\begin{array}{cc}
1 & 0 \\
0 & -1 
\end{array}\right) 
\quad \textit{and} \quad
\RAE=\left(\begin{array}{cc}
\frac{\rt3}{2} & \frac{1}{2} \\
\frac{1}{2} &  -\frac{\rt3}{2}
\end{array}\right).
\end{equation}

\begin{figure}[htbp]
\begin{center}
\includegraphics{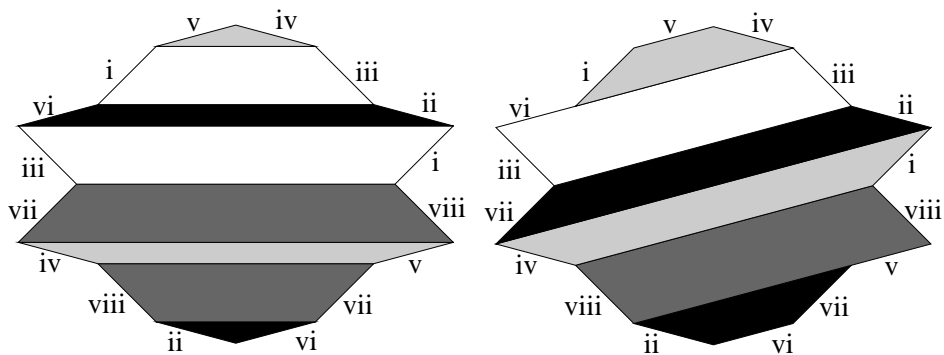}
\caption[Two decompositions into cylinders]{Decompositions into cylinders using saddle connections
parallel to lines $\ell_1$ and $\ell_2$ of figure \ref{fig:surface}.}
\label{fig:decomp12}
\end{center}
\end{figure}

Now we will find a parabolic automorphism of $\MT_\Delta$ fixing the point
$B$. There is a decomposition of $\MT_\Delta$ by saddle connections parallel
to $\ell_1$ of figure \ref{fig:surface}. This decomposition is depicted in figure \ref{fig:decomp12},
and cuts the surface into 4 cylinders. 
It can be verified that these 
cylinders come in pairs with two possible moduli:
\begin{equation}
\label{modB}
\frac{1}{5+3 \rt3} \quad \textit{and} \quad \frac{1}{10 + 6 \rt3}.
\end{equation}
In particular, note the first modulus is twice the second. Thus, there is a 
parabolic element of the automorphism group which fixes the horizontal
direction and acts as a single Dehn twist on the pair of 
cylinders with modulus $\frac{1}{10 + 6 \rt3}$ and a double Dehn twist on
the pair of cylinders with modulus $\frac{1}{5+3 \rt3}$. This parabolic is:
\begin{equation}
P_B=\left(\begin{array}{cc}
1 &  10 + 6 \rt3 \\
0 &  1
\end{array}\right). 
\end{equation}

Both the reflection $\RAB$ and the parabolic $P_B$ fix the point $B$,
so their composition does as well. Their composition gives another 
reflection:
\begin{equation}
\RBC=P_B \circ \RAB = \left(
\begin{array}{cc}
1 & -10-6\rt3 \\
0 & -1 
\end{array}
\right).
\end{equation}

The same idea will work for generating a reflection preserving points
$E$ and $D$. We decompose $\MT_\Delta$ into 4 cylinders using segments 
parallel to $\ell_2$ of figure \ref{fig:surface}. 
This decomposition is shown on the right in 
figure \ref{fig:decomp12}. The moduli of these cylinders again come in pairs: 
\begin{equation}
\label{modE}
\frac{3}{6+4\rt3} \quad \textit{and} \quad \frac{1}{6+4\rt3}.
\end{equation}
This means that there is a parabolic inducing a single Dehn twist 
on the cylinders with moduli $\frac{1}{6+4\rt3}$ and a triple Dehn twist
on the other cylinders. This parabolic is
\begin{equation}
P_E=\left(\begin{array}{cc}
\frac{1}{2}(-1-2\rt3) & \frac{1}{2}(12+7\rt3) \\
\frac{1}{2}(-\rt3) &  \frac{1}{2}(5+2\rt3)
\end{array}\right). 
\end{equation}
Again we obtain $\RDE$ by composition:
\begin{equation}
\RDE=\RAE \circ P_E = \left(\begin{array}{cc}
\frac{1}{2}(-3-\rt3) & \frac{1}{2}(13+7\rt3) \\
\frac{1}{2}(1-\rt3) &  \frac{1}{2}(3+\rt3)
\end{array}\right).
\end{equation}

\begin{figure}[htbp]
\begin{center}
\includegraphics{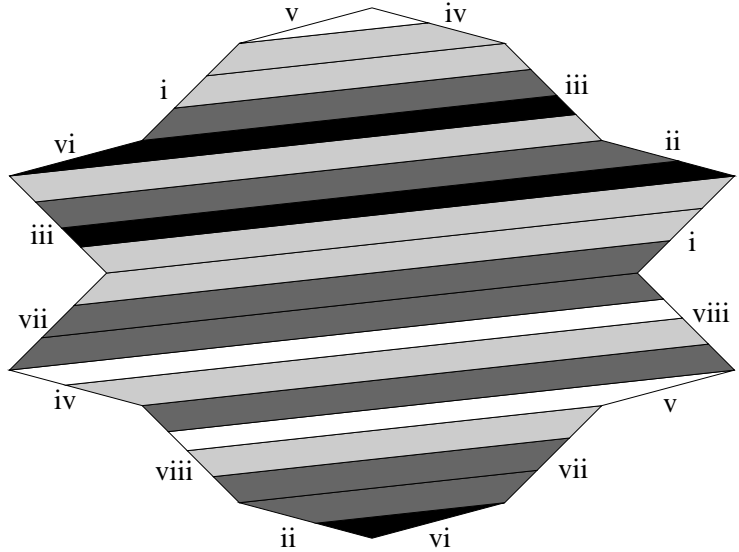}
\caption[Another decomposition into cylinders]{Saddle connections with slope $\frac{1}{11}(-4 + 3 \rt3)$
decompose the surface into these cylinders.}
\label{fig:decomp3}
\end{center}
\end{figure}

We apply the same trick one last time. $\RDE$ preserves
two parallel families of lines on the surface, each corresponding
to eigenvectors of the matrix. Of course, one is the family of
lines parallel to $\ell_2$. The second family has slope 
$\frac{1}{11}(-4 + 3 \rt3)$. We decompose the surface using saddle 
connections parallel to this direction (see figure \ref{fig:decomp3}). 
Again, these cut the surface into
four cylinders whose moduli come in pairs. The moduli are
\begin{equation}
\label{modD}
\frac{1}{29+17\rt3} \quad \textit{and} \quad \frac{1}{58+34\rt3}.
\end{equation}
Therefore, we get a parabolic fixing lines of slope 
$\frac{1}{11}(-4 + 3 \rt3)$:
\begin{equation}
P_D=\left(\begin{array}{cc}
\frac{1}{2}(-11-7\rt3) & \frac{1}{2}(115+67\rt3) \\
\frac{1}{2}(-1-\rt3) &  \frac{1}{2}(15+7\rt3)
\end{array}\right).
\end{equation}

We compose the parabolic with the reflection $\RDE$:
\begin{equation}
\RCD=\RDE \circ P_D = \left(\begin{array}{cc}
5+3 \rt3 & -51-30\rt3 \\
1 &  -5-3\rt3
\end{array}\right).
\end{equation}

One way to check that lines $\overline{BC}$ and $\overline{CD}$ 
intersect at a right angle is to compute the trace of the product of reflections
in the sides. If this trace is zero, then the sides meet at a right angle.
We compute
\begin{equation}
\textit{Tr}(\RBC \circ \RCD)=\textit{Tr}\left(\begin{array}{cc}
-5-3 \rt3 & 53+30\rt3 \\
-1 &  5+3\rt3
\end{array}\right) =0.
\end{equation}
$\diamondsuit$

\begin{lemma}
The reflections in the side of the polygon of figure \ref{fig:domain} together with $-I$ generate
the affine automorphism group $\Gamma(\MT_\Delta)$.
\end{lemma}
{\em Proof:} 
Let $G$ be the group generated by $-I$ and the reflections in the sides of the polygon of figure \ref{fig:domain}, and
let $G^+$ be the index 2 subgroup which preserves the orientation of $\H^2$.

Let $\Y_G=\H^2/G^+$, a sphere with 3 punctures and two cone singularities. One singularity has cone angle $\pi$ and the other has 
cone angle $\pi/3$. 
We can compute the area of $\Y_G$ using the Gauss-Bonnet formula. 
Recall the Gauss-Bonnet formula for hyperbolic surfaces with
cone singularities tells us that for a surface $S$ of genus $g$ with $p$ punctures and cone singularities of cone angle 
$\theta_1, \ldots, \theta_n$,
\begin{equation}
\textit{area}(S)=2 \pi (2g+p-2)+ \displaystyle\sum_{i=1}^n (2\pi-\theta_i)
\end{equation}
We compute that $\textit{area}(\Y_G)=\frac{14 \pi}{3}$.

Now let $V=\Gamma(\MT_\Delta)$ be the complete affine automorphism group, 
$V^+$ be the orientation preserving subgroup, and $\Y_V=\H^2/V^+$.
We wish to show $V=G$. The previous lemma showed that $G$ is a subgroup of $V$. 
Thus we have a covering map $\psi:\Y_G \to \Y_V$. Further we know that
\begin{equation}
[V^+:G^+]=\textit{area}(\Y_G)/\textit{area}(\Y_V)
\end{equation}
where $[V^+:G^+]$ is the index of the subgroup $G^+$ inside $V^+$. In order to show that
$V^+=G^+$, it is sufficient to show $\textit{area}(\Y_G)/\textit{area}(\Y_V) < 2$.

We would like to use Gauss-Bonnet on $\Y_V$. First we will show that $\Y_V$ also has $3$ punctures.
It is sufficient to show that none of the punctures of $\Y_G$ can be identified by $\psi$.
We will give affine invariants which distinguish the three decompositions into cylinders mentioned in the previous proof.
The ratio of the moduli of the cylinders associated to vertex $E$ of the polygons is $3$, while the ratios of the moduli
of cylinders associated to $B$ and $D$ are both $2$ (see equations \ref{modE}, \ref{modB}, and \ref{modD}). Thus
$E$ can not be identified with $B$ or $D$. Another affine invariant is the ratio of the widths of the cylinders.
We can compute that these ratios are
\begin{equation}
\begin{array}{ccc}
w_B=1+\sqrt{3} & \textrm{ and } & w_D=\frac{1+\sqrt{3}}{2}
\end{array}
\end{equation}
Thus, the punctures coming from $B$ and $D$ cannot be identified by the covering $\psi$. This shows that
$\Y_V$ has $3$ punctures.
 
We also need to show that $\Y_V$ has at least one cone singularity. The image of a cone singularity in
$\Y_G$ must be a cone singularity in $\Y_V$. Further, the image of the cone singularity with cone angle 
$\pi/3$ must be a cone singularity with cone angle $\theta$ which is less than $\pi/3$. Gauss-Bonnet now tells us
that 
\begin{equation}
\textit{area}(\Y_V) \ge 2 \pi (3-2) + (2 \pi - \theta) \ge 4 \pi - \pi/3 = 11 \pi / 3 > \frac{1}{2}  \textit{area}(\Y_G)
\end{equation}
Thus $\textit{area}(\Y_G)/\textit{area}(\Y_V)<2$, so $[V^+:G^+]=1$ and $V^+=G^+$. 

Finally, because both $V$ and $G$ contain orientation reversing elements, we know $[V:G]=[V^+:G^+]$. 
Thus $V=G$.
\endproof

\section{The other symmetric translation surface}
The translation surface $\MT_\Delta$ is affinely equivalent to another translation surface $\MT'$. This surface corresponds to the vertex $C$ of the fundamental domain in $\H^2$ of the affine automorphism group of $\MT_\Delta$ pictured in figure \ref{fig:domain}. Thus, $\MT'$ can be written as $A(\MT_\Delta)$ where $A \in \SLhat(2, \R)$ is
$$A=\left(
\begin{array}{cc}
1 & 5+3 \rt3 \\
0 & 1 
\end{array}
\right)$$

\begin{figure}[h]
\begin{center}
\includegraphics{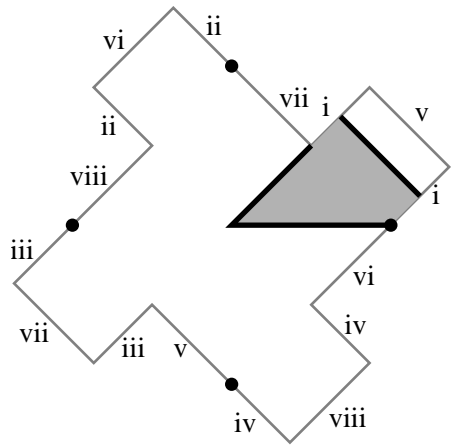}
\caption[Figure \ref{fig:surface}]{The translation surface $\MT_\Delta$ is affinely equivalent to this translation surface,
$\MT'$. The grey region is an annulus which arises as the quotient of this surface by its Euclidean isometries.}
\label{fig:other_table}
\end{center}
\end{figure}

The surface $\MT'$ is pictured in figure \ref{fig:other_table} and supports a Euclidean isometry group $G$, which is a 
dihedral group of order $8$. The quotient $\MT'/G$ is an annulus, which we will now describe. 
Consider the points
$$
\begin{array}{c}
P_0=(0,0) \quad P_1=(1,0) \quad P_2=(\frac{3+\rt3}{4},\frac{-1+\rt3}{4}) \\
P_3=(\frac{1+\rt3}{4},\frac{1+\rt3}{4}) \quad P_4=(\frac{1}{2},\frac{1}{2}).
\end{array} $$
$\MT'/G$ is the pentagon $P_0 P_1 P_2 P_3 P_4$ with edge $\overline{P_1 P_2}$ glued to $\overline{P_4 P_3}$ by translation.

\vspace{1em}

\noindent
{\bf Acknowledgements.}
I would like to thank Barak Weiss for the idea of searching for obtuse
Veech triangles, and Curt McMullen for his encouragement to write
this result down. I would also like to thank Rich Schwartz for his collaboration on McBilliards and my 
graduate advisor Yair Minsky.

\bibliography{/home/pat/active/my_papers/bibliography}
\end{document}